\documentclass[10pt,a4paper]{article}

\pdfoutput=1

\usepackage[utf8]{inputenc}
\usepackage[T2A]{fontenc}
\usepackage[bulgarian,english]{babel}
\usepackage{amssymb}
\usepackage{graphicx}
\usepackage{hyperref}

\newenvironment{items}[1][\labelitemi]{\begin{list}{#1}{\setlength{\topsep}{0pt}\setlength{\partopsep}{0pt}\setlength{\parsep}{0pt}\setlength{\itemsep}{\parskip}}}{\end{list}}

\def\fr#1#2{\leavevmode\kern.1em\raise.76ex\hbox{\scalebox{.9}{#1}}\kern-.15em/\kern-.1em\lower.40ex\hbox{\scalebox{.9}{#2}}}

\def\tpc#1{\vskip.7ex\noindent\textbf{\S#1.}\texttt{~}}

\title{\Large\sc Fraction Space Revisited\kern.5mm%
}

\author{\sc Boyko B.\ Bantchev}

\date{}

\begin{document}

\maketitle

\thispagestyle{empty}

{\narrower\small\noindent
Rationals are known to form interesting and computationally rich structures, such as Farey sequences and infinite trees.  Little attention is being paid to more general, systematic exposition of the basic properties of fractions as a set.  Some concepts are being introduced without motivation, some proofs are unnecessarily artificial, and almost invariably both seem to be understood as related to specific structures rather than to the set of fractions in general.  Surprisingly, there are essential propositions whose very statement seem to be missing in the number theory literature.  This article aims at improving on the said state of affairs by proposing a general and properly ordered exposition of concepts and statements about them.  In addition, historical remarks are made on generating the set of all fractions -- a much older discovery than it is widely believed.

}

\section*{Introduction}

The discovery of the Farey sequences in the beginning of the $19^\mathrm{th}$ century marked a shift of interest with regard to rational numbers from purely numeric properties to ones concerning the structure of the set of rationals.  It became apparent that not only could we add, multiply, compare, etc.\ fractions, but we can relate them to each other, including generate them from one another, based on new concepts -- ones of adjacency and computation of mediants.  Later investigations lead to the emergence of a construct now known as the Stern-Brocot (SB) tree, which is more general by including all simple fractions at once, and enables deeper exploration of the structure of $\mathbb{Q}$ than it is possible with Farey sequences.  Later still, even more structural relations of generative nature have been found in $\mathbb{Q}$.

Remarkably, the Farey sequences and the SB tree, although of different origins, exploit the same relational and computational basis -- adjacency and mediants.  Also, and somewhat surprisingly, the several other trees that have been discovered -- independently and in different contexts -- to contain the positive rationals, all appear to be indistinguishable as sets of rows.  The author sees this as an evidence that a more general and comprehensive theory of the structure of the set of rationals is awaiting its development.

As a small step in this direction, the author finds it necessary to emancipate the treatment of some fundamental concepts from the consideration of Farey sequences, SB tree or any similar structure.  Indeed, introduction of adjacency and mediants, along with the related properties, is at present tied to the said structures \cite{hardy,gkp}, leaving the false impression that it belongs there.  In fact, it need not be so.  Treating the basic concepts and their properties independently uncovers their true, fundamental meaning and fosters a more complete account of them.

The remainder of the text contains two parts.  The first is a brief exposition of terms and properties of simple fractions.  Introducing normalized error and distance first, mediants come as a solution to a naturally arising problem, and adjacency is the minimality of distance.  Computation of mediants of adjacent fractions is then shown to generate all positive rationals.  The distances between a fraction and a pair of adjacent fractions encode how the former is generated from the latter.  As a way of contrast, the reader is invited to compare this approach to a `Farey-based' text, such as \cite{hardy}, Ch.\,III or similar, or to \cite{dick}, end of Ch.V, or to \cite{gkp}~4.5, where both Farey sequences and the Stern-Brocot tree are discussed.

Facts pertaining to Farey sequences, the SB tree or other specific structures either receive a more general treatment here, or are not discussed but can easily be inferred from the rest.

The second part contains remarks on the history and results of generating and enumerating fractions, especially in a tree form.  Besides drawing the reader's attention to several important and useful enumerations, my goal is to stress that, contrary to popular belief, that history is rooted in the antiquity.

\section*{The space of fractions}

To save space, many proofs are omitted or only sketched in the following, especially if they are straightforward enough to reconstruct on one's own.

Following \cite{gkp}, the co-primality of integers $p$ and $q$ is denoted by $p\,\bot\,q$.

We consider fractions $\fr{$p$}{$q$}$ with $p,q\,{>}\,0$ and $p\,\bot\,q$.

\tpc{1}%
Fractions are commonly used to approximate numbers.  A number $x$ is approximated by a fraction $\fr{$p$}{$q$}$ with an absolute error $|x\,{-}\,\fr{$p$}{$q$}|$.  Instead of the latter, a \emph{normalized error}, $q\,|x\,{-}\,\fr{$p$}{$q$}|=|q\,x\,{-}\,p|\,,$ may be preferred for the following reasons.  The best absolute error is within \fr1{$2q$}; by multiplying with $q$, the impact of the denominator is neutralized -- as if integers are considered in place of fractions.  Also, a larger $q$ is more `resource consuming' (greater number of fractions to consider, with larger numerators and denominators), so it makes sense to `tax' for that by assuming a proportionally larger error.

\tpc{2}%
For similar reasons, for fractions $\fr{$p$}{$q$}\,{<}\fr{$r$}{$s$}$, in place of the absolute difference $\,\fr{$r$}{$s$}\,{-}\,\fr{$p$}{$q$}=\fr{$qr{-}ps$}{$qs$}\,$, by multiplying with $qs$, we introduce a normalized difference, or a \emph{distance} $|\fr{$p$}{$q$}\,,\!\fr{$r$}{$s$}|=qr{-}ps$.  The distance is $\ge$1, and when it is 1, we say that the fractions are \emph{adjacent} and write $\fr{$p$}{$q$}\,\bot\fr{$r$}{$s$}$\,.

It can immediately be observed that:
\begin{items}[\textbf{--}]
\item the notion of adjacency agrees with, and generalizes the one of integers (successive integers are adjacent fractions);
\item for any two fractions \fr{$p$}{$q$} and \fr{$r$}{$s$}\,:  \,\,$q\,r\,{-}\,p\,s\,=1\,\,\Leftrightarrow\,\,\fr{$p$}{$q$}\,{<}\fr{$r$}{$s$}\,\land\,\fr{$p$}{$q$}\,\bot\fr{$r$}{$s$}\,$;
\item $\fr{$p$}{$q$}\,\bot\fr{$r$}{$s$}\,\,\Rightarrow\,\,p\,\bot\,q\,\land\,p\,\bot\,r\,\land\,r\,\bot\,s\,\land\,q\,\bot\,s$;
\item $|\fr{$p$}{$q$}\,,\fr{$r$}{$s$}|=|\fr{$s$}{$r$}\,,\fr{$q$}{$p$}|=|\fr{$p$}{$r$}\,,\fr{$q$}{$s$}|=|\fr{$s$}{$q$}\,,\fr{$r$}{$p$}|$\,.
\end{items}

\tpc{3}%
Keeping in mind the second observation above, what if in $q\,r{-}p\,s\,{=}1$ any of $p$, $q$, $r$, and $s$ can be 0 and not only positive?

\textbf{Proposition}.
If $p,q,r,s\,{\ge}\,0$ and $q\,r{-}p\,s\,{=}1$, then exactly one of the following is true (as can be proved by case analysis and use of $1\,{=}\,qr{-}ps=q(r{-}p)+p(q{-}s)$):
\begin{items}
\item[0.] $p\,{=}\,s\,{=}\,0\,\land\,q\,{=}\,r\,{=}\,1$;
\item[1.] $p\,{=}\,r\,{=}\,1\,\land\,q\,{=}\,s{+}1$;
\item[2.] $q\,{=}\,s\,{=}\,1\,\land\,r\,{=}\,p{+}1$;
\item[3.] either \,$p\,{<}\,r\,\land\,q\,{<}\,s$\, or \,$p\,{>}\,r\,\land\,q\,{>}\,s$.
\end{items}

Thus we arrive at a useful generalization of adjacency: from now on, we consider $\,q\,r{-}p\,s\,{=}1\,$ to be the definition of adjacency of \fr{$p$}{$q$} and \fr{$r$}{$s$} even when numerators or denominators are 0.  In fact, only two new fractions, \fr01 and \fr10, become legitimate -- according to the above, they are the only admissible representatives of 0 and $\infty$, respectively (no \fr0{$n$} or \fr{$n$}0 for $n\,{>}1$ can be adjacent with any other fraction).

More generally, the sign of $q\,r{-}p\,s$ defines whether $\fr{$p$}{$q$}\,{<}\,\fr{$r$}{$s$}$\, or vice versa, and $q\,r{-}p\,s$ is, as before, the distance $|\fr{$p$}{$q$}\,,\!\fr{$r$}{$s$}|$ between fractions -- now including \fr01 and \fr10\,.

The 0.$-$2.\ cases from the proposition can be rewritten:
\begin{items}
\item[0.] $\fr{$p$}{$q$}\,{=}\fr01\,\land\,\fr{$r$}{$s$}\,{=}\fr10$;
\item[1.] $\fr{$p$}{$q$}\,{=}\fr1{$n$+1}\,\land\,\fr{$r$}{$s$}\,{=}\fr1{$n$}$ for $n\ge0$;
\item[2.] $\fr{$p$}{$q$}\,{=}\fr{$n$}1\,\land\,\fr{$r$}{$s$}\,{=}\fr{$n$+1}1$ for $n\ge0$.
\end{items}

Adjacency of fractions where at least one of them is an integer can be summarized as follows.
\begin{items}[\textbf{--}]
\item An integer $p$ is only adjacent with fractions of the kind $p\,{-}\fr1{$n$}$ from the left, and with \fr10 and fractions of the kind $p\,{+}\fr1{$n$}$ from the right.
In particular, $p$ is adjacent with $p{-}1$ and $p{+}1$, and 0 is adjacent with all fractions \fr1{$n$}, including \fr10, and only them.
\item \fr10 is adjacent with all fractions \fr{$n$}1 (the integers) and only them.
\end{items}

\tpc{4}%
Let $x\in[\fr{$p$}{$q$}\,,\fr{$r$}{$s$}]$.  The larger the $x$, the larger the (normalized) error \,$q\,x{-}p$\, of approximating $x$ with \fr{$p$}{$q$}, and the smaller the error \,$r{-}s\,x$\, of approximating $x$ with \fr{$r$}{$s$}.  The former error starts from 0, while the latter ends at 0.  The unique $M\in[\fr{$p$}{$q$}\,,\fr{$r$}{$s$}]$ for which the errors are equal is the solution for $x$ of the equation $\,q\,x{-}p\,{=}\,r{-}s\,x\,$ -- a fraction.  The same $M$ is also the unique solution for a fraction $f$ of the equation $\,|\fr{$p$}{$q$},f|\,{=}\,|f,\fr{$r$}{$s$}|$.

Therefore, the following holds:

\textbf{Proposition}.
For any fraction-bound interval $[\fr{$p$}{$q$}\,,\fr{$r$}{$s$}]$, the fraction $\,M{=}\,\fr{$p$+$r$}{$q$+$s$}\,$ is the unique one that is at equal distances from the two ends of the interval, and also the unique number that is approximated with equal errors by them.

It is immediately seen that the distances satisfy $\,|\fr{$p$}{$q$},M|\,{=}\,|M,\fr{$r$}{$s$}|\,{=}\,|\fr{$p$}{$q$}\,,\fr{$r$}{$s$}|\,$ and the error of approximation is $|\fr{$p$}{$q$}\,,\fr{$r$}{$s$}|/(q+s)$.

For any two fractions \fr{$p$}{$q$} and \fr{$r$}{$s$}, the fraction \fr{$p$+$r$}{$q$+$s$} is called their \emph{mediant}.
Due to the context in which the mediant emerged above, clearly $\fr{$p$}{$q$}\,{<}\,\fr{$p$+$r$}{$q$+$s$}\,{<}\,\fr{$r$}{$s$}$ holds.

\tpc{5}%
\textbf{Proposition}.
If $\fr{$p$}{$q$}\,\bot\fr{$r$}{$s$}$, then:
\begin{items}[\textbf{--}]
\item $p{+}r\,\bot\,q{+}s$ also holds, i.\,e.\ the mediant as computed is in lowest terms;
\item $\fr{$p$}{$q$}\,\bot\,\fr{$p$+$r$}{$q$+$s$}$\, and \,$\fr{$p$+$r$}{$q$+$s$}\,\bot\,\fr{$p$}{$q$}$ (and no other fraction is adjacent with both \fr{$p$}{$q$} and \fr{$r$}{$s$}).
\end{items}

\tpc{6}%
From \S3 we know that for any pair of adjacent fractions \fr{$p$}{$q$} and \fr{$r$}{$s$} except \fr01 and \fr10, there holds either $p\,{\le}\,r$ and $q\,{\le}\,s$ or $p\,{\ge}\,r$ and $q\,{\ge}\,s$.  Accordingly, we can construct either \fr{$r{-}p$}{$s{-}q$}\, or \,\fr{$p{-}r$}{$q{-}s$}, and name it \emph{medidifference}.  Like the mediant, the medidifference is adjacent to its formative fractions.  If $p\,{\le}\,r$ and $q\,{\le}\,s$, then \fr{$r$}{$s$} is the mediant of \fr{$p$}{$q$} and \fr{$r{-}p$}{$s{-}q$} (and similarly for the opposite inequalities).

\textbf{Proposition}.
If $p\,{<}\,r$ and $q\,{<}\,s$, then:
\begin{items}[\textbf{--}]
\item $\fr{$p$}{$q$}\,{<}\,\fr{$r$}{$s$}$\,,\,\,\,
$\fr{$p$}{$q$}\,{<}\,\fr{$p$+$r$}{$q$+$s$}$\,,\,\,\,
$\fr{$p$+$r$}{$q$+$s$}\,{<}\,\fr{$r$}{$s$}$\,,\,\,\,
$\fr{$p$}{$q$}\,{<}\,\fr{$r$$-$$p$}{$s$$-$$q$}$\,,\,\,\,
$\fr{$r$}{$s$}\,{<}\,\fr{$r$$-$$p$}{$s$$-$$q$}$\,,\,\,\,
and \\
\,\,\,$\fr{$p$+$r$}{$q$+$s$}\,{<}\,\fr{$r$$-$$p$}{$s$$-$$q$}$\,\,\,
are all equivalent;
\item $\fr{$p$}{$q$}\,\bot\fr{$r$}{$s$}$,\,\,\,
$\fr{$p$}{$q$}\,\bot\,\fr{$p$+$r$}{$q$+$s$}$,\,\,\,
$\fr{$r$}{$s$}\,\bot\,\fr{$p$+$r$}{$q$+$s$}$,\,\,\,
$\fr{$p$}{$q$}\,\bot\fr{$r$$-$$p$}{$s$$-$$q$}$,\,\,\,
$\fr{$r$}{$s$}\,\bot\fr{$r$$-$$p$}{$s$$-$$q$}$\,,\,\,
and \\
\,\,\,$|\fr{$p$+$r$}{$q$+$s$},\fr{$r$$-$$p$}{$s$$-$$q$}|\,{=}\,2$\,\,\,
are all equivalent.
\end{items}

\tpc{7}%
The statement of \S5 suggests that each interval whose bounds are adjacent fractions can be indefinitely subdivided by means of constructing mediants: the interval is first split into two, then into four, etc., each time obtaining new adjacent pairs whose corresponding intervals can be further subdivided.  Each new mediant is different from all others, as it belongs strictly within an interval which, by the moment the mediant is computed, does not contain other fractions.

\tpc{8}%
Just as mediant computation applies for subdividing an interval, the reverse process takes place by computing medidifferences.  Based on the results from \S6, an interval bounded by adjacent fractions can be extended to the left or right, then the same can be done with the extended interval, etc.\ -- as many times as needed, or until $[\fr01\,,\fr10]$ is reached.  Thus e.\,g.\ from $[\fr43\,,\fr75]$ we obtain $[\fr43\,,\fr32]$, $[\fr11\,,\fr32]$, $[\fr11\,,\fr21]$, $[\fr11\,,\fr10]$, and $[\fr01\,,\fr10]$.

In essence, Euclid's algorithm for finding GCD is being simultaneously run on two pairs of natural numbers.  Each pair is co-prime (see \S2), and the simultaneity is possible due to \S6.

\tpc{9}%
\textbf{Proposition}.
For each fraction $\fr{$a$}{$b$}\in[\fr{$p$}{$q$}\,,\fr{$r$}{$s$}]$ there holds
\,$da\,{=}\,np{+}mr$, \,$db\,{=}\,nq{+}ms$,\, where 
\,$d\,{=}\,|\fr{$p$}{$q$}\,,\fr{$r$}{$s$}|$,\, $m\,{=}\,|\fr{$p$}{$q$}\,,\fr{$a$}{$b$}|$,\, and \,$n\,{=}\,|\fr{$a$}{$b$}\,,\fr{$r$}{$s$}|$.
(Can be verified by immediate checking.)

The above means, in particular, that $\fr{$a$}{$b$}\,{=}\fr{$np{+}mr$}{$nq{+}ms$}$\, and that when $\fr{$p$}{$q$}\,\bot\fr{$r$}{$s$}$\,, the right-hand side is in lowest terms.

\textbf{Proposition}.
If $\fr{$p$}{$q$}\,\bot\fr{$r$}{$s$}$\, and $m$ and $n$ are as above, then $m\,\bot\,n$.

Indeed, as $\fr{$p$}{$q$}$ and $\fr{$r$}{$s$}$ are adjacent, then so are $\fr{$p$}{$r$}$ and $\fr{$q$}{$s$}$\,.  Applying Euclid's algorithm simultaneously to the pairs $p,q$ and $r,s$ as in \S8, then multiplying the results by $m$ and $n$, respectively, and adding them up, we would be thus applying the same algorithm to the pair $a,b$.  At some point, $p,q$ would reduce to $0,1$ and $r,s$ to $1,0$.  The respective values to which $a$ and $b$ reduce would then be $m$ and $n$.  As $a\,\bot\,b$, there must also be $m\,\bot\,n$.

The integers $m$ and $n$, besides being distances, can also be regarded as \emph{coordinates} of \fr{$a$}{$b$} with respect to $[\fr{$p$}{$q$}\,,\fr{$r$}{$s$}]$.

The following statements can be proved by immediate verification.

\textbf{Proposition}.
The mediant of two fractions $f_1$ and $f_2$ within $[\fr{$p$}{$q$}\,,\fr{$r$}{$s$}]$ has a `coordinate fraction' \fr{$m$}{$n$} which is itself the mediant of the coordinate fractions of $f_1$ and $f_2$ (in lowest terms).

\textbf{Proposition}.
If \,$\fr{$a$}{$b$}\,{=}\fr{$n_1p{+}m_1r$}{$n_1q{+}m_1s$}$\,\, and \,$\fr{$c$}{$d$}\,{=}\fr{$n_2p{+}m_2r$}{$n_2q{+}m_2s$}$\,,\, then \\
\,$|\fr{$a$}{$b$}\,,\fr{$c$}{$d$}\,|\,{=}\,|\fr{$p$}{$q$}\,,\fr{$r$}{$s$}|.|\fr{$m_1$}{$n_1$},\fr{$m_2$}{$n_2$}|$.

\noindent
If $\fr{$p$}{$q$}\,\bot\fr{$r$}{$s$}$\,, then the above simplifies to \,$|\fr{$a$}{$b$}\,,\fr{$c$}{$d$}\,|\,{=}\,|\fr{$m_1$}{$n_1$},\fr{$m_2$}{$n_2$}|$, i.\,e.\ the distance between fractions equals the one between their coordinates.

\tpc{10}%
\textbf{Proposition}.
The process of mediant subdivision of an interval $[\fr{$p$}{$q$}\,,\fr{$r$}{$s$}]$, bounded by adjacent fractions, generates all the fractions within the interval.

For any fraction $f\in[\fr{$p$}{$q$}\,,\fr{$r$}{$s$}]$, take $m$ and $n$ as in \S9.  If $m\,{\le}\,n$, then $f\,{=}\,((n-m)p{+}m(p{+}r))/((n-m)q{+}m(s{+}q))$.  Thus $f$, with coordinates $m,n$, is obtainable by subdividing $[\fr{$p$}{$q$}\,,\fr{$r$}{$s$}]$ exactly when it, with coordinates $m,n{-}m$ in the narrower interval $[\fr{$p$}{$q$}\,,\fr{$p${+}$r$}{$q${+}$s$}]$, is obtainable by subdividing that interval (its bounds also being adjacent).  Similarly when $m\,{\ge}\,n$.

Repeating the same reduction, we are once more applying Euclid's algorithm, here to the pair $m,n$, while transforming the way $f$ is being represented.  Eventually, one of the two coordinates is zeroed, while in the currently reached interval $f$ becomes trivially representable, as it coincides with one of its ends.

It follows that, starting from $[\fr01,\fr10)$, all nonnegative rationals can be obtained by subdividing with mediants.  Every integer $\fr{$n$}1>0$ is generated as the mediant of an integer (the preceding one) and \fr10\,:\ from \fr01 and \fr10 we get \fr11, from \fr11 and \fr10 -- \fr21, and so on.

The process of fraction generation through mediant computation repeats itself in many respects and forms, such as the following:

$\cdot$ For any integer $n$, fraction generation within $[n,n{+}1]$ repeats that within $[0,1]$.

$\cdot$ Due to the symmetry between numerators and denominators of mediants, fraction generation within $[\fr11,\fr10)$ and that within $(\fr01,\fr11]$ reciprocally mirror each other.

$\cdot$ For any interval $[\fr{$p$}{$q$}\,,\fr{$r$}{$s$}]$ with $\fr{$p$}{$q$}\,\bot\fr{$r$}{$s$}$, in view of \fr{$p$}{$q$} and \fr{$r$}{$s$} having coordinate fractions \fr01 and \fr10 with respect to the interval, the third proposition of \S9 establishes a duality between the fractions within that interval and those within $[\fr01,\fr10)\equiv[0,\infty)$.  In other words, there is a natural 1-1 correspondence between the fractions within any specific adjacent-bound interval and the set of all nonnegative fractions.

\tpc{11}%
Adjacent fractions can be found as closely as needed:

\textbf{Proposition}.
For each fraction $\fr{$p$}{$q$}\in(\fr01\,,\fr10)$ and each $\varepsilon\,{>}\,0$, there exist \fr{$a$}{$b$} and \fr{$c$}{$d$}, such that
$\fr{$a$}{$b$}\,{<}\,\fr{$p$}{$q$}\,{<}\,\fr{$c$}{$d$}$\,,
$\fr{$p$}{$q$}\,{-}\fr{$a$}{$b$}\,{<}\,\varepsilon$,
$\fr{$c$}{$d$}\,{-}\fr{$p$}{$q$}\,{<}\,\varepsilon$,
$\fr{$a$}{$b$}\,\bot\fr{$p$}{$q$}$\,,
and
$\fr{$p$}{$q$}\,\bot\fr{$c$}{$d$}$\,.

\tpc{12}%
\textbf{Proposition}.
The coordinate representation of fractions with respect to an interval with adjacent ends $[\fr{$p$}{$q$}\,,\fr{$r$}{$s$}]$ entails that the denominators and the numerators of such fractions are no smaller than $q{+}s$ and $p{+}r$, respectively.  In fact, they are strictly greater everywhere except at the mediant.

A reverse statement also holds.

\textbf{Proposition}.
If every fraction within $(\fr{$p$}{$q$}\,,\fr{$r$}{$s$})$ (assuming $\fr{$r$}{$s$}\,{<}\fr10$) has a denominator \,${>}\max(q,s)$ (or, equivalently, has a numerator \,${>}\max(p,r)$), then $\fr{$p$}{$q$}\,\bot\fr{$r$}{$s$}$ holds.

(Take $\fr{$a$}{$b$}\,{=}\fr01$ and $\fr{$c$}{$d$}\,{=}\fr10$ and narrow $[\fr{$a$}{$b$}\,,\fr{$c$}{$d$}]$ by mediant subdividing until either (a)~$\fr{$a$}{$b$}\,{=}\,\fr{$p$}{$q$}$ and $\fr{$r$}{$s$}\,{=}\,\fr{$c$}{$d$}$ -- so indeed $\fr{$p$}{$q$}\,\bot\fr{$r$}{$s$}$ -- or (b)~one of \fr{$a$}{$b$} and \fr{$c$}{$d$} is strictly in, and the other strictly out of $(\fr{$p$}{$q$}\,,\fr{$r$}{$s$})$.  In case of (b), let e.\,g.\ $\fr{$a$}{$b$}\,{<}\,\fr{$p$}{$q$}\,{<}\,\fr{$c$}{$d$}\,{<}\,\fr{$r$}{$s$}$.  Then, due to $\fr{$a$}{$b$}\,\bot\fr{$c$}{$d$}$\,, $q\,{>}\,d$ holds, which is a contradiction with \fr{$c$}{$d$} being within $(\fr{$p$}{$q$}\,,\fr{$r$}{$s$})$. Thus (b) is not possible.)

The above two propositions can be united in the following (again, assuming $\fr{$r$}{$s$}\,{<}\fr10$)

\textbf{Proposition}.
$\fr{$p$}{$q$}\,\bot\fr{$r$}{$s$}\,\,\Leftrightarrow\,\,\forall\,\fr{$a$}{$b$}\in[\fr{$p$}{$q$}\,,\fr{$r$}{$s$}]:\,b\,{>}\max(q,s)$.

\tpc{13}%
\textbf{Proposition}.
Every positive fraction is the mediant of exactly one pair of nonnegative adjacent fractions.
(From \S10 we know that at least one pair exists; assuming that there are more than one can be shown to produce a contradiction.)

Regarding the fractions whose mediant is a given fraction as its parents, the proposition says that every positive fraction has unique parents, which can be distinguished as a left and right one.  From \S6 it follows that one of these parents is itself a parent of the other, thus being a parent and a grandparent of the given fraction.  A fraction can be regarded as `left' or `right' depending on which of its parents is also its grandparent.

The uniqueness of parents has many consequences.  For example, any two adjacent fractions belong to a closed interval $[n,n{+}1]$, where $n\,{\ge}\,0$ is an integer, including when one or both fractions are precisely $n$ or/and $n{+}1$.  (The opposite would mean that an integer exists between the fractions, and can therefore be obtained through subdividing the respective interval; this contradicts to an integer being a mediant of another integer and \fr10.)  Also, the parents of a fraction that is not an integer belong to the same interval between two successive integers as the fraction itself.  (The parents belong to one such interval, and that is necessarily the same that the mediant belongs to.)

The above proposition establishes a fundamental fact of fractions.  Yet, to the best of the author's knowledge, it -- along with a number of others listed here -- is left out from the literature on elementary number theory.  The closest similar theorems concern the uniqueness with respect to appearance in Farey sequences or something else, but these are weaker propositions.  The omission can be attributed to the generally insufficient thoroughness that the subject has been treated with thus far.

\section*{Trees and enumerations: historical remarks}

Farey sequences were the first and for some time the only observed manifestation of the structure inherent to the set of fractions.  Their significance lessened with the discovery of several binary trees, each of which contains the whole set of positive fractions.  Such a tree is attractive by providing a single, uniform structure in which all fractions are represented and none is repeated.  Consequently, the entire set of fractions can be enumerated, mapped to other sets, etc.  The figure shows the first several rows of some trees: each fraction is considered the parent (in the tree) of the two immediately below it.

The first all-fractions tree that became known is the Stern-Brocot tree, named after the mathematician M.\,Stern and the watchmaker A.\,Brocot.  M.\,Stern investigated `diatomic sequences' of integers, previously (1850) defined in a more general setting by G.\,Eisenstein \cite{stern}.  Those sequences were further studied by others; an account of the obtained results can be found in \cite{giuli}.  Neither Stern nor the others considered fractions, however.

Contrastingly, A.\,Brocot used a process which was effectively mediant insertion, to approximate any given fraction with fractions of lesser denominators \cite{broc}.  This inspired several French mathematicians from the end of the $19^\mathrm{th}$ century to study sequences of fractions, similar to Farey's, but somewhat more general.

Of course, none of Stern, Brocot, or the other authors mentioned trees explicitly, trees being a concept introduced in mathematics only several decades ago now, and closely related to computing and algorithmics.  Nevertheless, Stern's sequences clearly exhibit hierarchical structure.  The SB tree itself seems to have been introduced by \cite{gkp}.  On the figure, each fraction is a mediant of the fractions above it that are its closest to the left and to the right; e.\,g., \fr47 is the mediant of \fr12 and \fr35.  (It is assumed that \fr01 and \fr10 are above all other fractions, in a leftmost and rightmost positions, respectively.)

A similar tree became known as the Calkin-Wilf (CW) tree, after the authors of \cite{calk}.  In it, a fraction \fr{$p$}{$q$} has a left child \fr{$p$}{$p$+$q$} and a right child \fr{$p$+$q$}{$q$}\,.  Calkin and Wilf also pointed out that the tree, i.\,e.\ all positive rationals, can be enumerated by the function $f(1)=1, \; f(2\,n)=f(n), \; f(2\,n+1)=f(n)+f(n+1)$ in the sense that the ratios $f(n)/f(n+1)$ for consecutive $n$ are in reduced form and run through all tree nodes row by row, each row in a left to right order.

Interestingly, the same function was studied earlier in \cite{dijk}, but in a different context and without its enumeration property.  More importantly, the function has been in use since even earlier (see e.\,g.\ \cite{lind} and \cite{rezn}), as a generator of M.\,Stern's sequence of integers!  In this respect, \cite{back} points out that in Stern's work the CW tree is even more readily recognized than what is now called the SB tree, and consequently, the former tree should more appropriately be named after Eisenstein and Stern.

The recurrence in the definition of the above function means that enumerating the rationals using it requires an amount of memory that increases with the argument.  However, a formula was discovered that directly, in constant space, generates all fractions in the same order as $f$: if $x$ is any fraction, its successor in the sequence is $\,x'=(\lfloor{}x\rfloor+1-\{x\})^{-1}=(2\lfloor{}x\rfloor+1-x)^{-1}$.  (According to \cite{knuth}, the author of this result is M.\,Newman.)

Yet another fraction tree is defined in \cite{shen} and independently in \cite{andr}.  (Neither actually speaks of trees, but the hierarchy is apparent from the defining formula.)  A fraction \fr{$p$}{$q$} in this tree has a left child \fr{$p$+$q$}{$q$} and a right child \fr{$q$}{$p$+$q$}\,.  It can be called Shen-Andreev, or SA tree.

The three mentioned trees are constructed according to very different principles, so they may seem to have very little in common.  In one respect, though, they are remarkably similar: they are `row-wise equivalent', i.\,e.\ they only differ in how the fractions are permuted \emph{within rows}.  In particular, each fraction is on the same row with its reciprocal -- this is obvious for the SB tree due to its symmetric structure, and therefore is true for all the three.  Thus, on every row, there are as many rationals $<$\kern.07em{}1 as there are $>$1.

One may observe that for enumeration purposes it may suffice to consider only rationals less than 1.  Indeed, pairing each such number with its reciprocal readily turns any enumeration of these numbers into an enumeration of all rationals.  If from any of the above three trees all fractions $>$1 and \fr11 itself are removed (and the remaining positions are considered adjacent), we obtain a binary tree rooted at \fr12 and containing all positive fractions $<$\kern.07em{}1 and only them.  From the SB tree we would thus obtain its left sub-tree, and from the SA tree -- one usually called `the Kepler tree', where it can be seen that \fr{$p$}{$q$} has a left child \fr{$p$}{$p$+$q$} and a right child \fr{$q$}{$p$+$q$}\,.  The result from similarly reducing the CW tree does not seem to have appeared in the literature.

Kepler's name was given to the tree because the latter appears in book~{\small III}, chapter~{\small II} of his magnum opus \cite{kepler}, in relation with computing the harmonic divisions of a string.  Kepler actually constructed a finite, small portion of the complete tree, and did not comment on whether that tree would be exhaustive and nonrepetitive of the fractions $<$\kern.07em{}1.  He nevertheless gave the construction rule and was clear about evolving \fr12 into a binary hierarchy.

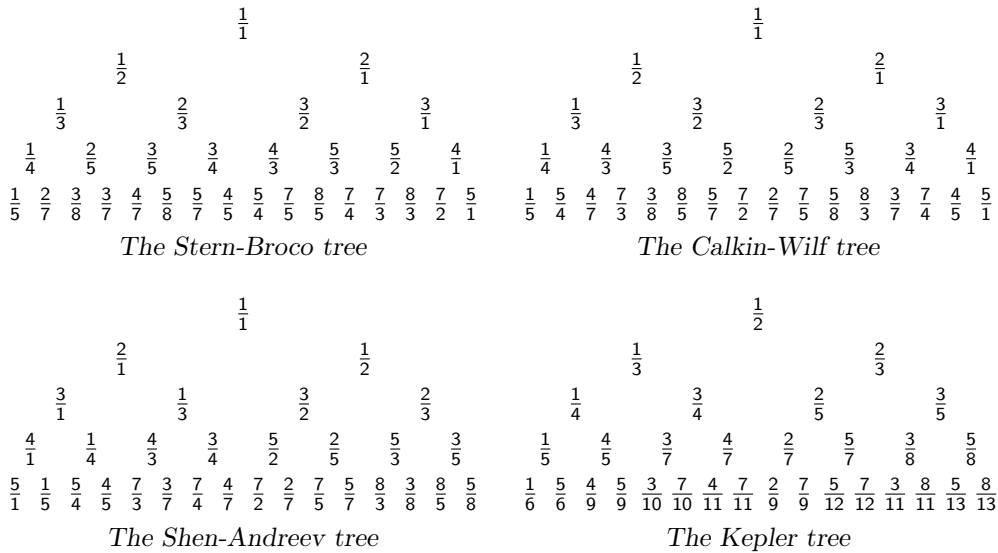
\def\rr#1#2#3{\hbox to #3{\hfil$\mathsf{\frac{#1}{#2}}$\hfil}}
\newdimen\sbtbw \sbtbw=2.0mm
\def\vssbt{\vskip.5em}
\def\riv#1#2{\rr{#1}{#2}{2\sbtbw}}
\def\riii#1#2{\rr{#1}{#2}{4\sbtbw}}
\def\rii#1#2{\rr{#1}{#2}{8\sbtbw}}
\def\ri#1#2{\rr{#1}{#2}{16\sbtbw}}
\def\rz#1#2{\rr{#1}{#2}{32\sbtbw}}
\def\osbt{\vrule width 0pt\kern1.4mm}
\begin{center}
\begin{minipage}[t]{33.3\sbtbw}
{\parindent=0pt
\vbox{%
\noindent
\osbt\rz{1}{1} \vssbt
\osbt\ri{1}{2}\ri{2}{1} \vssbt
\osbt\rii{1}{3}\rii{2}{3}\rii{3}{2}\rii{3}{1} \vssbt
\osbt\riii{1}{4}\riii{2}{5}\riii{3}{5}\riii{3}{4}\riii{4}{3}\riii{5}{3}\riii{5}{2}\riii{4}{1} \vssbt
\osbt\riv{1}{5}\riv{2}{7}\riv{3}{8}\riv{3}{7}\riv{4}{7}\riv{5}{8}\riv{5}{7}\riv{4}{5}\riv{5}{4}\riv{7}{5}\riv{8}{5}\riv{7}{4}\riv{7}{3}\riv{8}{3}\riv{7}{2}\riv{5}{1}
}}
\vspace{-.5em}\begin{center}\slshape The Stern-Broco tree\end{center}
\end{minipage}
\begin{minipage}[t]{33.3\sbtbw}
{\parindent=0pt
\vbox{%
\noindent
\osbt\rz{1}{1} \vssbt
\osbt\ri{1}{2}\ri{2}{1} \vssbt
\osbt\rii{1}{3}\rii{3}{2}\rii{2}{3}\rii{3}{1} \vssbt
\osbt\riii{1}{4}\riii{4}{3}\riii{3}{5}\riii{5}{2}\riii{2}{5}\riii{5}{3}\riii{3}{4}\riii{4}{1} \vssbt
\osbt\riv{1}{5}\riv{5}{4}\riv{4}{7}\riv{7}{3}\riv{3}{8}\riv{8}{5}\riv{5}{7}\riv{7}{2}\riv{2}{7}\riv{7}{5}\riv{5}{8}\riv{8}{3}\riv{3}{7}\riv{7}{4}\riv{4}{5}\riv{5}{1}
}}
\vspace{-.5em}\begin{center}\slshape The Calkin-Wilf tree\end{center}
\end{minipage}
\vskip1.5em
\begin{minipage}[t]{33.3\sbtbw}
{\parindent=0pt
\vbox{%
\noindent
\osbt\rz{1}{1} \vssbt
\osbt\ri{2}{1}\ri{1}{2} \vssbt
\osbt\rii{3}{1}\rii{1}{3}\rii{3}{2}\rii{2}{3} \vssbt
\osbt\riii{4}{1}\riii{1}{4}\riii{4}{3}\riii{3}{4}\riii{5}{2}\riii{2}{5}\riii{5}{3}\riii{3}{5} \vssbt
\osbt\riv{5}{1}\riv{1}{5}\riv{5}{4}\riv{4}{5}\riv{7}{3}\riv{3}{7}\riv{7}{4}\riv{4}{7}\riv{7}{2}\riv{2}{7}\riv{7}{5}\riv{5}{7}\riv{8}{3}\riv{3}{8}\riv{8}{5}\riv{5}{8}
}}
\vspace{-.5em}\begin{center}\slshape The Shen-Andreev tree\end{center}
\end{minipage}
\begin{minipage}[t]{33.3\sbtbw}
{\parindent=0pt
\vbox{%
\noindent
\osbt\rz{1}{2} \vssbt
\osbt\ri{1}{3}\ri{2}{3} \vssbt
\osbt\rii{1}{4}\rii{3}{4}\rii{2}{5}\rii{3}{5} \vssbt
\osbt\riii{1}{5}\riii{4}{5}\riii{3}{7}\riii{4}{7}\riii{2}{7}\riii{5}{7}\riii{3}{8}\riii{5}{8} \vssbt
\osbt\riv{1}{6}\riv{5}{6}\riv{4}{9}\riv{5}{9}\riv{3}{10}\riv{7}{10}\riv{4}{11}\riv{7}{11}\riv{2}{9}\riv{7}{9}\riv{5}{12}\riv{7}{12}\riv{3}{11}\riv{8}{11}\riv{5}{13}\riv{8}{13}
}}
\vspace{-.5em}\begin{center}\slshape The Kepler tree\end{center}
\end{minipage}
\end{center}

Discovering organizational patterns in the set of the rational numbers, and the related enumeration(s), seem too important to be neglected.  Yet, historically, neglecting has been the rule rather than an exception.  That the discovery of such simple objects as mediants and Farey sequences took so very long is already strange enough.  And would they ever be widely known if it were not for the influence of Hardy\,\&\,Wright's classic \cite{hardy}, written well more than a century later?

Also, even today the enumerability of $\mathbb{Q}$ is almost ubiquitously demonstrated by referring to the Cantors's method – despite Stern's work predating Cantor's (1895) and being more relevant to the case.  Each of the SB etc.\ trees establishes an explicit one-to-one correspondence between $\mathbb{N}$ and $\mathbb{Q}$.  On the contrary, Cantor's proof (\cite{cantor}, p.\,107) constructs such a correspondence between $\mathbb{N}$ and $\mathbb{N\times{}N}$, from which there only follows the \textsl{existence}, and not an actual construction, of a one-to-one mapping between $\mathbb{N}$ and $\mathbb{Q}$.

An even more striking evidence of the undeserved lack of interest the mathematicians have had regarding the structure of rationals is the fact that a construct equivalent to the CW tree has been in existence since at least as early as the $3^\mathrm{rd}$ century~\textsc{bce}.  It has been before the eyes of the mathematical world and yet apparently remained unnoticed!

Nicomachus of Gerasa \cite{nicom} and Theon of Smyrna \cite{theon}, both works written in 1-2 century~\textsc{ce}, describe a process for generating all proportions starting from equality.  Theon says that it has been known to Eratosthenes of Cyrene, so its author must have lived perhaps in the $4^\mathrm{th}$ century~\textsc{bce}.  The computation starts from $\,1\;\;1\;\;1\,$ and, in general, from $\,x\;\;y\;\;z\,$ it obtains $\,x\;\;x{+}y\;\;x{+}2\,y{+}z\,$ and $\,x{+}2\,y{+}z\;\;y{+}z\;\;z\,$, each of which it then transforms in the same way, etc.

As \cite{shch} observes, each triple thus obtained is of the form $\,a^2\;\;ab\;\;b^2\,$, from where it is obvious that it represents the same ratio $\,a^2\,{:}\,ab=ab\,{:}\,b^2=a\,{:}\,b\,$ redundantly, and therefore one can replace triples with pairs, so that the above process amounts to transforming each pair $\,a\;\,b$, starting from $\,1\;\;1$, into $\,a\;\,a{+}b\,$ and $\,a{+}b\;\,b$, resulting in a binary hierarchy of ratios $\,a\,{:}\,b$.  Of course the latter is nothing else but the CW tree (although the author of \cite{shch}, being more a philosopher than a mathematician, does not mention this).  It is clear now that the generation of the rational numbers through a simple process, and the CW tree in particular, are not only not new but a valuable achievement of the ancient Greek mathematics.

Both Nicomachus's and Theon's works are well known.  Specifically, the former has been translated to Latin by Boethius and published as \cite{boet} (6${}^\mathrm{th}$ century), serving as the standard work on arithmetic throughout the middle ages, and remaining influential up to the end of the 18${}^\mathrm{th}$ century.  Although briefly, the proportion generation process is also described in another widely known and respected work -- Pappus of Alexandria's (4${}^\mathrm{th}$ century) \cite{papp} -- book~{\small III}, prop.~{\small XVII \& XVIII}.  In view of all this, it is amazing for how long the construction of the tree of rationals has escaped the attention of number theorists.

Historians of mathematics have also been silent on the subject.  Of Heath, Struik, Neugebauer, Dickson, van der Waerden, and several others, there is only one to mention it at all, and the only comment he has is: `None of this is very profound, but it is rather nice' \cite{waer}.  He, too, seems to have not fully understood the meaning behind the text that he discussed.

\end{document}